%% file: main.tex
\documentclass{article}
\usepackage{graphicx} 
\usepackage{prodint}
\usepackage{amsmath,amssymb,amsthm}
\usepackage{xcolor}
\usepackage[round]{natbib}
\usepackage[scale=0.75]{geometry} 

\theoremstyle{plain}
\newtheorem{defi}{Definition}
\newtheorem{theorem}{Theorem}
\newtheorem{lemma}{Lemma}
\newtheorem{remark}{Remark}
\newtheorem{example}{Example}
\newtheorem{cor}{Corollary}

\newcommand{\E}{\mathrm{E}}

\title{Conditional Delta-Method \\ for Resampling Empirical Processes \\ in Multiple Sample Problems}
\author{Merle Munko\footnote{Department of Mathematics, Otto-von-Guericke University Magdeburg}, \ \
Dennis Dobler\footnote{Department of Statistics, TU Dortmund University; \newline \phantom{!} \quad Research Center Trustworthy Data Science and Security, University Alliance Ruhr.}}
\date{\today}

\begin{document}

\maketitle

\begin{abstract}
    The functional delta-method has a wide range of applications in statistics. 
    Applications on functionals of empirical processes yield various limit results for classical statistics. 
    To improve the finite sample properties of statistical inference procedures that are based on the limit results, resampling procedures such as random permutation and bootstrap methods are a popular solution.
    In order to analyze the behaviour of the functionals of the resampling empirical processes, corresponding conditional functional delta-methods are desirable.
    While conditional functional delta-methods for some special cases already exist, there is a lack of more general conditional functional delta-methods for resampling procedures for empirical processes, such as the permutation and pooled bootstrap method.
    This gap is addressed in the present paper. 
    Thereby, a general multiple sample problem is considered.
    The flexible application of the developed conditional delta-method is shown in various relevant examples.
\end{abstract}

\input{manuscript}

\section*{Acknowledgements}
We want to thank Eric Beutner for helpful discussions.
Furthermore, Merle Munko gratefully acknowledges support from the \textit{Deutsche Forschungsgemeinschaft} (grant no. DI 2906/1--2 and GRK 2297 MathCoRe). 

\bibliographystyle{plainnat} 
\bibliography{main}

\end{document}

%% file: manuscript.tex
\section{Introduction}

Many applications of statistics involve comparisons of multiple samples. 
Section 3.8 of the monograph by \cite{vdVW23} is devoted to a related empirical process treatment.
In addition to an analysis of differences of two independent empirical processes, they also explained how to analyze a random permutation and a pooled bootstrap version of the empirical processes.
Most statistical applications involve a functional that is applied to these empirical processes.
The statistical properties of a functional of one empirical process can be derived with the help of the functional delta-method; cf.\ Section~3.10 in \cite{vdVW23}.
The application to a multiple sample problem, however, requires a delta-method with a varying reference point.

Several extensions of the functional delta-method in different directions have already been investigated in the literature.
For example, \cite{fang2019inference} studied the inference of functionals that are only directionally differentiable
and, recently, \cite{neumeyer2023generalizedhadamarddifferentiabilitycopula} proposed a generalization of Hadamard differentiability for applications of the functional delta-method to the empirical copula processes.
Under measurability assumptions, \cite{BEUTNER20102452} developed a modified functional delta method for quasi-Hadamard differentiable functionals.
Additionally, there exist conditional delta-methods for the bootstrap (in one sample) in Section~3.10.3 of \cite{vdVW23} and also extensions on (uniformly) quasi-Hadamard differentiable functionals for the bootstrap under measurability assumptions \citep{BeutnerZaehleBootstrap, beutner2016functionaldeltamethodbootstrapuniformly}.
However, as far as we know, a two- or multiple sample equivalent of such delta-methods for resampling empirical processes is not available in the literature. In detail, most of the existing methods require some of the following:
\begin{enumerate}
    \item[(i)] measurability assumptions,
    \item[(ii)] that the resampling counterpart converges weakly to the same limit as the empirical process, and 
    \item[(iii)] a fixed centering element of the empirical process, particularly independent of the sample sizes.
\end{enumerate}
However, these requirements are usually not satisfied for resampling methods for empirical processes in multiple sample problems, such as random permutation and pooled bootstrapping.

In this paper, we will develop a conditional delta-method in outer probability without assuming (i)--(iii) for applications to the randomly permuted and pooled bootstapped empirical processes in multiple independent sample problems. 
To this end, we require the uniform Hadamard differentiability of the functionals applied to the empirical processes.
In several examples, we show its applicability and usefulness. This includes conditional central limit theorems for the permutation and pooled bootstrap counterparts of the Wilcoxon statistic, the Nelson-Aalen estimator and the Kaplan-Meier estimator.

The remainder of this paper is organized as follows. 
Section~\ref{sec:Model} includes four subsections: 
in Section~\ref{ssec:Notation}, the model of our multiple sample problem is presented and the notation used in this paper is introduced. 
Existing convergence results of the resampling empirical processes are restated in Section~\ref{ssec:Weak}.
Additionally, a limit theorem for the permutation empirical process of multiple samples as an extension of Theorem~3.8.1 in \cite{vdVW23} is developed.
A functional delta-method for the empirical processes of the multiple sample problem is obtained in Section~\ref{ssec:Functional}.
In Section~\ref{ssec:Uniform}, uniform Hadamard differentiability is defined and some properties are investigated. 
Section~\ref{sec:Results} contains the main results of this paper that cover a flexible conditional delta-method. 
Particularly, this delta-method is applicable for deriving the limit of functionals of permutation and pooled bootstrap empirical processes. 
Exemplary functionals, applications, and limitations of our main result are given in Section~\ref{sec:Examples}. 
This includes the Wilcoxon functional in Section~\ref{ssec:Wilcoxon}, the product integral in Section~\ref{ssec:Product}, and the inverse map in Section~\ref{ssec:Inverse}.
We conclude with a summary and a discussion of possibilities for future research in Section~\ref{sec:Discussion}.

\section{Theoretical Framework and Preliminaries}\label{sec:Model}

\subsection{Model and Notation}\label{ssec:Notation}
Let, for each of $m \geq 2$ independent samples, $\boldsymbol{X}_{11}, \dots, \boldsymbol{X}_{1n_1} \sim P_1$, $\dots, $ $\boldsymbol{X}_{m1}, \dots, \boldsymbol{X}_{mn_m} \sim P_m$, be i.i.d.\ random elements on a measurable space $(\chi, \mathcal{A})$ and let $\mathbb{P}_{j, n_j} := \frac 1 {n_j} \sum_{i=1}^{n_j} \delta_{\boldsymbol{X}_{ji}} $ be the $j$-th empirical measure, $j=1, \dots, m$.

The introduction of the resampling techniques for the empirical process requires the pooled data. 
To this end, denote the pooled sample by 
$$(\boldsymbol{Z}_{N1}, \dots, \boldsymbol{Z}_{NN}) := (\boldsymbol{X}_{11}, \dots, \boldsymbol{X}_{1 n_1}, \dots, \boldsymbol{X}_{m1}, \dots, \boldsymbol{X}_{mn_m}), $$
where $N:= \sum_{j=1}^m n_j$ is the total sample size.
Let $\boldsymbol{R} = (R_1, \dots, R_N)$ be a vector that is uniformly distributed on the set of all permutations of $\{1,2, \dots, N\}$ and independent of the data $\boldsymbol{Z}_{N1}, \dots, \boldsymbol{Z}_{NN}$. 
Also, let $N_j := \sum_{\ell=1}^j n_\ell$ be the total sample size of the first $j$ samples, with $N_0 := 0$.
Then, the multiple sample permutation empirical measures are defined as
$$ \mathbb{P}_{j,n_j}^\pi := \frac1 {n_j}\sum_{i=N_{j-1}+1}^{N_j} \delta_{\boldsymbol{Z}_{N R_i}}, \quad j=1, \dots, m. $$

Next, denote by $\mathbb{H}_N := \frac1N \sum_{i=1}^N \delta_{\boldsymbol{Z}_{Ni}} = \sum_{j=1}^m \frac{n_j}N \mathbb P_{j,n_j}$ the pooled empirical measure.
The multiple sample bootstrap empirical measures
are defined as
$$ \hat{\mathbb{P}}_{j,n_j} := \frac1 {n_j}\sum_{i=N_{j-1}+1}^{N_j} \delta_{\hat{\boldsymbol{Z}}_{N i}} , \quad j=1, \dots, m,$$
where $\hat{\boldsymbol{Z}}_{N 1},...,\hat{\boldsymbol{Z}}_{N N} \sim \mathbb{H}_N$ is an i.i.d.\ sample drawn from the pooled empirical measure.


Throughout the paper, we assume that $\frac{n_j}N \to \lambda_j \in (0,1)$, as $\min_{j=1,\dots, m} n_j \to\infty$. 
Denote $\boldsymbol{\lambda} := (\lambda_1, \dots, \lambda_m)$ and $H := \sum_{i=1}^m \lambda_i P_i$.
Furthermore, let $\mathcal{F}$ denote a class of measurable functions $f: \chi \to \mathbb{R}$ that is $P_i$-Donsker for all $i=1,...,m$, i.e., $\sqrt{n_i}(\mathbb{P}_{i,n_i} - P_i) \rightsquigarrow \mathbb{G}_i$ in the space $\ell^{\infty}(\mathcal{F})$ of all bounded real-valued functions on $\mathcal F$ as $n_i \to\infty$, where here and throughout $\mathbb{G}_i$ is a tight $P_i$-Brownian bridge for all $i=1,...,m$ and $\rightsquigarrow$ denotes weak convergence in the sense of Section~1.3 in \cite{vdVW23}.
In the following, let $\mathbb{G}_1,\dots,\mathbb{G}_m$ be independent.
An immediate consequence of the above is
\begin{align}
    \label{eq:Hn}
    \sqrt{N} (\mathbb H_N - H_{\boldsymbol{n}}) =  \sum_{j=1}^m\frac{\sqrt{n_j}}{\sqrt{N}}\sqrt{n_j} (\mathbb P_{j,n_j} - P_j)\rightsquigarrow \mathbb G_{\boldsymbol{\lambda}} \quad\text{in $\ell^{\infty}(\mathcal{F})$}
\end{align}
as $\min_{j=1, \dots, m} n_j \to \infty$, where
 $H_{\boldsymbol{n}} := \sum_{j=1}^m \frac{n_j}N P_j$ and $\mathbb G_{\boldsymbol{\lambda}} := \sum_{j=1}^m \sqrt{\lambda_j} \mathbb{G}_j $.
 It should be noted that the centering element $H_{\boldsymbol{n}}$ in the previous display generally depends on the sample sizes.

\subsection{Weak Convergence Results of Resampling Empirical Processes}\label{ssec:Weak}
Now, we turn to the asymptotic behaviour of the resampling empirical processes.
We will see that the limits of the permutation and pooled bootstrap empirical process generally do not coincide with the limit of the empirical processes or the pooled empirical process.

In the two-sample case $m=2$, Theorem~3.8.1 in \cite{vdVW23} yields under $||P_i||_{\mathcal{F}} < \infty, i = 1,2,$ that $\sqrt{n_1} (\mathbb{P}^{\pi}_{1,n_1} - \mathbb{H}_N) \rightsquigarrow \sqrt{1-\lambda_1} \mathbb{G}_H$ given $\boldsymbol{X}_{11}, \boldsymbol{X}_{12}, \dots, $ $\boldsymbol{X}_{21}, \boldsymbol{X}_{22}, \dots$ in outer probability, where $\mathbb{G}_H$ denotes a tight $H$-Brownian bridge on $\ell^{\infty}(\mathcal{F})$. 
Here and throughout, stated convergence results are always meant as $\min_{j=1,\dots, m} n_j \to\infty$.
In Lemma~S.6 in the supplement of \cite{ditzhaus2021qanova}, the almost sure version of this theorem is generalized for multiple samples. Here, we state the corresponding extension in probability which is sufficient for most statistical applications.

\begin{theorem}\label{Permutation}
    Let $\mathcal{F}$ satisfy $||P_i||_{\mathcal{F}} < \infty$ for all $i=1,...,m$. Then, we have $\sqrt{N} (\mathbb{P}^{\pi}_{1,n_1} - \mathbb{H}_N,\dots,\mathbb{P}^{\pi}_{m,n_m} - \mathbb{H}_N) \rightsquigarrow \mathbb{G}_H^{\pi}$ given
    \begin{align}\label{eq:data}
        \boldsymbol{X}_{11}, \boldsymbol{X}_{12}, \dots, \boldsymbol{X}_{21}, \boldsymbol{X}_{22}, \dots, \dots, \boldsymbol{X}_{m1}, \boldsymbol{X}_{m2}, \dots
    \end{align}
    in outer probability, where $\mathbb{G}_H^{\pi}$ denotes a tight zero-mean Gaussian process on $(\ell^{\infty}(\mathcal{F}))^m$ with covariance function 
    $\boldsymbol\Sigma_H^{\pi}: {\mathcal{F}}^{m \times m} 
    \to \mathbb{R}^{m\times m} $. The component functions of $\boldsymbol\Sigma_H^{\pi}$ at $(f,g) = ((f_1,\dots,f_m),(g_1,\dots,g_m))$ are given by
    $$\left(\boldsymbol\Sigma_H^{\pi}(f,g)\right)_{ij} := \left(\lambda_i^{-1}{1}\{i=j\} -1 \right)H(f_i-Hf_i)(g_j-Hg_j)$$
    for all $i,j=1,\dots,m$, where here and throughout $1\{...\}$ denotes the indicator function of the set $\{...\}$.
\end{theorem}
\begin{proof}
   By imitating the proof of Theorem~3.8.1 in \cite{vdVW23}, we obtain conditional weak convergence    $\sqrt{N} (\mathbb{P}^{\pi}_{i,n_i} - \mathbb{H}_N) \rightsquigarrow \sqrt{(1-\lambda_i)/\lambda_i} \mathbb{G}_H$ given \eqref{eq:data} in outer probability for all $i=1,\dots,m$.  Then, the conditional asymptotic equicontinuity in probability of the vector $\sqrt{N} (\mathbb{P}^{\pi}_{1,n_1} - \mathbb{H}_N,\dots,\mathbb{P}^{\pi}_{m,n_m} - \mathbb{H}_N)$ follows easily by the conditional asymptotic equicontinuity in probability of the components $\sqrt{N} (\mathbb{P}^{\pi}_{i,n_i} - \mathbb{H}_N), i=1,\dots,m$. Furthermore, $\mathcal{F}$ is totally bounded in $L_2(P_i)$ for all $i=1,\dots,m$ by the assumptions.  
   
   Hence, it remains to consider the marginal distributions. Therefore, we proceed similar to the proof of (S.18) in the supplement of \cite{ditzhaus2021qanova} with a Cramér-Wold argument, where $P = H$ and $\mathbb{P} = \mathbb{H}_N$.
   Let $g_i = c_i f_i$ with $c_i\in [-1,1] $ and $f_i \in \mathcal{F}$ for all $i=1,\dots,m$. 
   Then, $$ \frac{1}{n_i} \max\{g_r(\boldsymbol X_{ij})^2 : j=1,\dots,n_i\} \to 0$$ holds almost surely for all $r,i=1,\dots,m$, which can be shown with the three steps (i)--(iii) in the beginning of the proof of Lemma~S.6 in the supplement of \cite{ditzhaus2021qanova} by using $g_r$ instead of the envelope function $\Tilde{F}$. 
   In detail, the three steps are the following:  
   \begin{enumerate}
       \item[(i)] dividing $g_r$ into $g_{r,1,M} := g_r 1\{|g_r| \leq M\}$ and $g_{r,2,M} := g_r 1\{|g_r| > M\}$  for $M\in\mathbb{N}$,
       \item[(ii)] using the inequalities
$(a + b)^2 \leq 2a^2 + 2b^2$ and 
$$\max_{j=1,\dots,n_i} g_{r,2,M}(\boldsymbol X_{ij})^2 \leq \sum_{j=1}^{n_i} g_{r,2,M}(\boldsymbol X_{ij})^2,$$ 
       \item[(iii)] letting first $n\to\infty$ and finally
$M\to\infty$.
   \end{enumerate}
   In the last step (iii) we need $P_ig_r^2 < \infty$ for an application of the dominated convergence theorem as $M\to\infty$, which holds due to the assumption that $\mathcal{F}$ is $P_i$-Donsker. Moreover, we have $\mathbb{H}_Ng_r \to Hg_r$ almost surely for all $r=1,\dots,m$ by the strong law of large numbers since $\mathcal{F}$ is $P_i$-Donsker for all $i=1,\dots,m$, and, thus, $H$-Donsker. 
   Condition (S.19) in the supplement of \cite{ditzhaus2021qanova} follows almost surely by the same arguments as given there.
   For proving condition (S.20) almost surely, it remains to show that $\mathbb{H}_N(g_ig_r) \to H(g_ig_r)$ almost surely for all $i,r=1,\dots,m$ by the last display in the proof of Lemma~S.6 in the supplement of \cite{ditzhaus2021qanova}. Since $\mathcal{F}$ is $P_i$-Donsker for every $i=1, \dots, m$, 
   $H(g_ig_r)$ exists and, thus, the almost sure convergence follows by the strong law of large numbers.
   Hence, the marginal convergence follows from (S.18) in the supplement of \cite{ditzhaus2021qanova} given \eqref{eq:data} almost surely.
   
\end{proof}

For the bootstrap empirical measure, \cite{vdVW23} showed in Theorem~3.8.6 for the two-sample case that $||P_i||_{\mathcal{F}} < \infty, i = 1,2,$ implies $\sqrt{n_1}(\hat{\mathbb{P}}_{1,n_1}-\mathbb{H}_N) \rightsquigarrow \mathbb{G}_H$ given $\boldsymbol{X}_{11}, \boldsymbol{X}_{12}, \dots, \boldsymbol{X}_{21}, \boldsymbol{X}_{22}, \dots$ in outer probability. Due to the independence, it follows easily that $||P_i||_{\mathcal{F}} < \infty, i = 1,\dots,m,$ implies 
\begin{align}\label{eq:Bootstrap}
    \sqrt{N}(\hat{\mathbb{P}}_{1,n_1}-\mathbb{H}_N,\dots,\hat{\mathbb{P}}_{m,n_m}-\mathbb{H}_N) \rightsquigarrow (\lambda_1^{-1/2}\mathbb{G}_{H,1},\dots,\lambda_m^{-1/2}\mathbb{G}_{H,m})
\end{align}
given \eqref{eq:data} in outer probability, where $\mathbb{G}_{H,1},\dots,\mathbb{G}_{H,m}$ denote independent $H$-Brownian bridges.

\subsection{Functional Delta-Method in the Multiple Sample Problem}\label{ssec:Functional}
In statistical applications, usually a functional is applied to the empirical processes. 
Delta-methods can be used to analyze the asymptotic behaviour of the functionals of empirical processes.
A delta-method for the empirical processes of the multiple sample problem can be easily obtained by applying Theorem 3.10.4 of \cite{vdVW23}, but we wished to make it more explicit, tailored to the problem at hand. 

\begin{theorem}\label{thm2}
    Let $\mathbb E$ be a metrizable topological vector space and $\phi:(\ell^{\infty}(\mathcal F))^m \to \mathbb E$ 
    such that $$\sqrt{N} (\phi(\mathbf{P} + N^{-1/2} \mathbf{h}_{\boldsymbol{n}}) - \phi(\mathbf{P})) \to \phi'_{\mathbf{P}}(\mathbf h) $$
as $\min_{j=1, \dots, m} n_j \to \infty$ holds for every converging sequence $\mathbf h_{\boldsymbol{n}} \to \mathbf h \in (\ell^{\infty}(\mathcal F))^m$ 
 with $\boldsymbol{n} := (n_1,\dots,n_m)$, $\mathbf{P} := (P_1,\dots,P_m)$, $\mathbf{P}+ N^{-1/2} \mathbf{h}_{\boldsymbol{n}} \in \mathbb (\ell^{\infty}(\mathcal F))^m$ for all $\boldsymbol{n}$
 and for an arbitrary map $\phi'_{\mathbf{P}}:(\ell^{\infty}(\mathcal F))^m \to \mathbb E$. Then, we have
$$ \sqrt{N} (\phi(\mathbb P_{1,n_1},\dots,\mathbb P_{m,n_m}) - \phi(\mathbf{P}) ) \rightsquigarrow \phi'_{\mathbf{P}}(\lambda_1^{-1/2}\mathbb G_1,\dots, \lambda_m^{-1/2}\mathbb G_m). $$
 If $\phi'_{\mathbf{P}}$ is linear and continuous, the sequence $$\sqrt{N} (\phi(\mathbb P_{1,n_1},\dots,\mathbb P_{m,n_m}) - \phi(\mathbf{P}) ) - \phi'_{\mathbf{P}}( \sqrt{N} (\mathbb P_{1,n_1} - P_1,\dots, \mathbb P_{m,n_m} - P_m) )$$ converges to zero in outer probability.
\end{theorem}
\begin{proof}
    This is an application of Theorem 3.10.4 of \cite{vdVW23}.
\end{proof}

The condition on $\phi$ in the previous theorem is satisfied if $\phi$ is Hadamard differentiable at $\mathbf{P}$. To define Hadamard differentiability of a functional $\phi: \mathbb{D}_\phi \subset \mathbb D \to \mathbb E$, let $\mathbb D$ and $\mathbb E$ be metrizable topological vector spaces.

\begin{defi}[Hadamard differentiability]
    The functional $\phi$ is called Hadamard differentiable at $\theta\in\mathbb{D}_{\phi}$ tangentially to a subspace $\mathbb D_0 \subset \mathbb D$ if $$t_n^{-1} (\phi(\theta + t_n h_n) - \phi(\theta)) \to \phi'_{\theta}(h) $$
holds for all $t_n\to 0$ and every converging sequence $h_n \to h \in \mathbb{D}_0$
 with $\theta+ t_n h_n \in \mathbb D_\phi$ for all $n$ and for a continuous, linear map $\phi'_{\theta}:\mathbb{D}_0 \to \mathbb{E}$.
\end{defi}

In order to obtain a delta-method for the permutation and pooled bootstrap empirical processes, we need to introduce the uniform Hadamard differentiability in the following section.

\subsection{Uniform Hadamard differentiability}\label{ssec:Uniform}
We aim to develop functional delta-methods that are suitable for applications to 
$(\mathbb P_{1,n_1}^\pi, \dots, \mathbb P_{m,n_m}^\pi)$ and to $(\hat{\mathbb P}_{1,n_1}, \dots, \hat {\mathbb P}_{m,n_m})$, conditionally on \eqref{eq:data}.
To this end, we will consider again a functional
$\phi: \mathbb{D}_\phi \subset \mathbb D \to \mathbb E$, where $\mathbb D$ and $\mathbb E$ are metrizable topological vector spaces. 

\begin{defi}[Uniform Hadamard differentiability]
    The functional $\phi$ is called uniformly Hadamard differentiable at $\theta\in\mathbb{D}_{\phi}$ tangentially to a subspace $\mathbb D_0 \subset \mathbb D$ if $$t_n^{-1} (\phi(\theta_{n} + t_n h_{n}) - \phi(\theta_{n})) \to \phi'_{\theta}(h) $$
holds for all $t_n\to 0$ and every converging sequence $h_{n} \to h \in \mathbb{D}_0$ and $\theta_{n} \to \theta$
 with $\theta_{n}, \theta_{n}+ t_n h_{n} \in \mathbb D_\phi$ for all $n$ and for a continuous, linear map $\phi'_{\theta}:\mathbb{D}_0 \to \mathbb{E}$.
\end{defi}
If the subspace $\mathbb{D}_0$ is not specified, we assume $\mathbb{D}_0 = \mathbb{D}$ in the following.
For example, $\mathbb D$ can be chosen as product space $(\ell^\infty(\mathcal{F}))^m$ equipped with the max-sup-norm for applications on the empirical measures.

In the following, we investigate different properties of uniform Hadamard differentiable functionals.
The following remarks address the classical Hadamard derivative as a special case, the more restrictive Fr\'echet differentiability, and the aggregation of multiple functionals.

\begin{remark}\label{Remark1}
    Let $\phi: \mathbb{D}_\phi \subset \mathbb D \to \mathbb E$ be Hadamard differentiable at $\theta\in\mathbb{D}_{\phi}$ tangentially to a subspace $\mathbb D_0 \subset \mathbb D$ with Hadamard derivative $ \phi'_{\theta}:\mathbb{D}_0 \to \mathbb{E}$. If $\phi$ is uniformly Hadamard differentiable at $\theta\in\mathbb{D}_{\phi}$ tangentially to $\mathbb D_0$, then the (uniform) Hadamard derivative is $ \phi'_{\theta}$, which can easily be seen by setting $\theta_{n} = \theta$ in the definition.
\end{remark}

\begin{remark}[Uniform Fréchet differentiability and other sufficient criteria for uniform Hadamard differentiability]
    Let $(\mathbb D, ||.||_{\mathbb D})$ and $(\mathbb E, ||.||_{\mathbb E})$ be normed spaces.
    We call a functional $\phi$ 
    {uniformly Fr\'echet differentiable}
    at $\theta \in \mathbb D_\phi$ with continuous and linear derivative $\phi'_\theta$ if
    $$ \| \phi(\theta + h) - \phi(\theta + k) - \phi_\theta'(h - k)\|_{\mathbb E} = o(\|h-k\|_{\mathbb D}), \quad \text{as } \| h\|_{\mathbb D}, \|k\|_{\mathbb D} \to 0.$$
    To see that the uniform Fr\'echet differentiability implies the uniform Hadamard differentiability of $\phi$, insert $h =  \theta_{n}+ t_n h_{n} - \theta$, $k =  \theta_{n} - \theta$.
    According to a variant of Problem~3.10.1 in \cite{vdVW23}, uniform Fr\'echet differentiability at $\theta$ is implied by the Fr\'echet differentiability on a neighborhood of $\theta$ and the uniform norm-continuity of $\vartheta \mapsto \phi_\vartheta'$ at $\theta$ if there exists a convex neighborhood of $\theta$ as a subset of $\mathbb D_\phi$.
Other criteria for the uniform Hadamard differentiability of $\phi$ at $\theta$ are the convexity of $\mathbb D_\phi$, the Hadamard differentiability on a neighborhood of $\theta$, $\lim_{\eta \to \theta} \phi'_\eta(h) = \phi'_\theta(h)$ for every $h \in \textnormal{lin} \ \mathbb D_\phi$, and $\lim_{\eta \to \theta} \phi'_\eta(h) = \phi'_\theta(h)$ uniformly in $h \in K$ for every totally bounded subset $K \subset \mathbb D_\phi$, where $\eta \in \mathbb D_\phi$; cf.\ (3.10.6) in \cite{vdVW23}.
\end{remark}


\begin{remark}
    \label{vector}
    If $\phi_1:\mathbb{D}_\phi \to \mathbb{E}_1, \dots, \phi_m: \mathbb{D}_\phi \to \mathbb{E}_{m}$ are uniformly Hadamard differentiable at $\theta \in \mathbb{D}_\phi \subset \mathbb{D}$ tangentially to $\mathbb{D}_0 \subset \mathbb{D}$ with Hadamard derivatives $\phi_{1,\theta}':\mathbb{D}_0 \to \mathbb{E}_1,\dots,\phi_{m,\theta}':\mathbb{D}_0 \to \mathbb{E}_{m}$, it follows directly from the definition of uniform Hadamard differentiability that $\phi:= (\phi_1,\dots,\phi_m): \mathbb{D}_\phi \to \mathbb E_1 \times\dots\times \mathbb E_m$ is uniformly Hadamard differentiable at $\theta$ tangentially to $\mathbb{L}_0$ with Hadamard derivative $\phi'_{\theta} = (\phi_{1,\theta}',\dots,\phi_{m,\theta}'):\mathbb{D}_0 \to \mathbb{E}_1 \times\dots\times \mathbb E_m$.
    Here, the product space $\mathbb E_1 \times\dots\times \mathbb E_m$ is equipped with the product topology.
\end{remark}

Finally, the following theorem provides a chain rule for uniformly Hadamard differentiable functionals.
It should be noted that \cite{beutner2016functionaldeltamethodbootstrapuniformly} proved a chain rule for uniformly quasi-Hadamard differentiable functionals; see Lemma~A.1 therein.
That chain rule implies the chain rule statement below.
For the sake of completeness, however, we shall present a version of the chain rule which is relevant for the remainder of this paper.

\begin{theorem}[Chain rule]\label{chain}
Let $\mathbb{L}, \mathbb{D}, \mathbb{E}$ be metrizable topological vector spaces.
    If $\psi: \mathbb{L}_\psi \subset\mathbb{L} \to \mathbb{D}_\phi \subset \mathbb{D}$ is uniformly Hadamard differentiable at $\theta \in \mathbb{L}_\psi$ tangentially to $\mathbb{L}_0 \subset\mathbb{L}$ with Hadamard derivative $\psi'_{\theta}:\mathbb{L}_0\to\mathbb{D}$ and $\phi: \mathbb{D}_\phi \to \mathbb{E}$ is uniformly Hadamard differentiable at $\psi(\theta)$ tangentially to $\psi'(\mathbb{L}_0)$ with Hadamard derivative $\phi'_{\psi(\theta)}:\psi'(\mathbb{L}_0)\to\mathbb{E}$, then $\phi \circ \psi: \mathbb{L}_\psi \to \mathbb E $ is uniformly Hadamard differentiable at $\theta$ tangentially to $\mathbb{L}_0$ with derivative $\phi_{\psi(\theta)}' \circ \psi_\theta'$.
\end{theorem}
\begin{proof}
Let $t\to 0, h_t \to h\in\mathbb{L}_0, \theta_t\to\theta$ with $\theta_t,\theta_t+th_t \in \mathbb{L}_{\psi}$.
    Write 
    \begin{align*}
     &(\phi \circ \psi)(\theta_t + t h_t) - (\phi \circ \psi)(\theta_t)  = \phi (\psi(\theta_t) + t k_t) - \phi (\psi(\theta_t)) ,
    \end{align*}
    where $k_t = (\psi(\theta_t + t h_t) - \psi(\theta_t))/t$.
    Now, the uniform Hadamard differentiability of $\psi$ yields that $k_t \to \psi'_\theta(h)$, next the uniform Hadamard differentiabiliy implies the theorem.
\end{proof}

\section{Main Results}\label{sec:Results}
In this section, we aim to prove a conditional delta method, e.g., for applications to the permutation and pooled bootstrap empirical processes.
For proving the main theorem, we first need the following auxiliary lemma to obtain joint unconditional convergence of two maps.
The intuition is that $M_n$ will stand for additional randomness, e.g., induced by random permutation or pooled bootstrapping, whereas $X_n$ will represent the original data.
The result is similar to the results in Section~2 and~3 in \cite{bucher2019note} but allows arbitrary maps in general metric spaces.

\begin{lemma}\label{JointConvergenceLemma}
    Let $\mathbb D, \mathbb E$ be metric spaces, $X_n:\Omega_1\to\chi_{1n}, M_n:\Omega_2\to\chi_{2n}$ be sequences of functions, where $(\Omega_1\times\Omega_2, \mathcal{A}_1\times\mathcal{A}_2, Q_1 \otimes Q_2)$ denotes a product probability space, and $h_n:\chi_{1n} \to \mathbb D$ be such that $h_n(X_n) \rightsquigarrow H$ for some separable Borel measurable random element $H:\Omega_1 \to \mathbb D$. Moreover, let $y_n:\chi_{1n} \times \chi_{2n} \to \mathbb E$ with $y_n(X_n, M_n) \rightsquigarrow Y$ conditionally on $X_n$ in outer probability for some separable Borel measurable random element $Y:\Omega_1\times\Omega_2 \to \mathbb E$.
    Then, it follows that $(h_n(X_n), y_n(X_n,M_n)) \rightsquigarrow (H,Y)$ unconditionally for independent $H, Y$.
\end{lemma}

A precise formulation of $y_n(X_n, M_n) \rightsquigarrow Y$ conditionally on $X_n$ in outer probability is that
    \begin{align}\label{eq:precise}
    &\sup\limits_{h\in BL_1(\mathbb E)} \left| \E_2 h\left( y_n(X_n, M_n) \right)^{2*}  - \E h(Y) \right| \to 0,
       \\& \label{eq:precise2}
        \E_2 h\left( y_n(X_n, M_n) \right)^*  - \E_2 h\left( y_n(X_n, M_n) \right)_*  \to 0
    \end{align}
    in outer probability for all $h\in BL_1(\mathbb E)$,
    where $\E_2$ denotes the conditional expectation with respect to $M_n$ given $X_n$, $BL_M(\mathbb E)$ denotes the set of all real functions on $\mathbb E$ with a Lipschitz norm bounded by $M \geq 0$ and the asterisks denote the minimal measurable majorants and maximal measurable minorants with respect to $X_n, M_n$ jointly
    for $^*$ and with respect to $M_n$ for $^{2*}$.

\begin{proof}
    We aim to apply Corollary~1.4.5 in \cite{vdVW23}. Hence, it remains to show $\E^* f(h_n(X_n))g(y_n(X_n, M_n)) \to \E f(H) \E g(Y) $ for all bounded, nonnegative Lipschitz functions $f: \mathbb D \to \mathbb R, g: \mathbb E \to \mathbb R$.
    We have
    \begin{align*}
        &\left|\E^* f(h_n(X_n)) g(y_n(X_n, M_n)) - \E f(H) \E g(Y)\right|
        \\&\leq \left|\E_1\E_2 \left(f(h_n(X_n)) g(y_n(X_n, M_n))\right)^*  -   \E_1 \left(f(h_n(X_n))^*\E_2 g(y_n(X_n, M_n))^*\right)\right| \\& \quad + \left|\E_1 \left(f(h_n(X_n))^*\E_2 g(y_n(X_n, M_n))^*\right) - \E_1 f(H) \E g(Y)\right|       , 
    \end{align*} 
    where $\E_1,\E_2$ denote the expectations regarding $(\Omega_1,\mathcal A_1, Q_1)$ and $(\Omega_2,\mathcal A_2, Q_2)$, respectively. By Lemma~1.2.2~(v) in \cite{vdVW23} and the nonnegativity of $f$ and $g$, it holds that
    $$ \left(f(h_n(X_n)) g(y_n(X_n, M_n))\right)^* \leq f(h_n(X_n))^* g(y_n(X_n, M_n))^*$$ almost surely
and, hence, it follows that
$$ \E_2 \left(f(h_n(X_n)) g(y_n(X_n, M_n))\right)^* \leq f(h_n(X_n))^* \E_2 g(y_n(X_n, M_n))^*$$ almost surely.
Thus, we get
\begin{align}
        &\left|\E^* f(h_n(X_n)) g(y_n(X_n, M_n)) - \E f(H) \E g(Y)\right| \notag
        \\&\leq \E_1\left(\left(f(h_n(X_n))^*\E_2 g(y_n(X_n, M_n))^*\right) - \E_2 \left(f(h_n(X_n)) g(y_n(X_n, M_n))\right)^*\right) \notag  \\& \quad + \left|\E_1 \left(f(h_n(X_n))^*\left(\E_2 g(y_n(X_n, M_n))^* - \E g(Y)\right)\right)\right| \notag
        \\& \quad+ \left|\E_1\left(\left( f(h_n(X_n))^* -  f(H) \right)\E g(Y)\right)\right|     \notag
        \\&\leq \E_1\left(\left(f(h_n(X_n))^*\E_2 g(y_n(X_n, M_n))^*\right) -  f(h_n(X_n))_* \E_2g(y_n(X_n, M_n))_*\right) \label{eq:1} \\& \quad + ||f||_{\infty} \E_1
        \left( \E_2 g(y_n(X_n, M_n))^* - \E_2 g(y_n(X_n, M_n))^{2*}  + \left|\E_2 g(y_n(X_n, M_n))^{2*}  - \E g(Y)\right|\right)
        \label{eq:2}
        \\& \quad+ ||g||_{\infty} \left|\E_1\left( f(h_n(X_n))^* -  f(H) \right)\right|     \label{eq:3}        
        .
    \end{align} 
    By $g(y_n(X_n, M_n))_* \leq g(y_n(X_n, M_n))^{2*} \leq g(y_n(X_n, M_n))^*$, \eqref{eq:precise} and \eqref{eq:precise2} imply $$\E_2 g(y_n(X_n, M_n))^* - \E_2 g(y_n(X_n, M_n))^{2*}  + \left|\E_2 g(y_n(X_n, M_n))^{2*}  - \E g(Y)\right| \to 0$$ in outer probability. Hence, the dominated convergence theorem provides that \eqref{eq:2} converges to zero.
    Due to $h_n(X_n) \rightsquigarrow H$, \eqref{eq:3} converges to zero. Hence, \eqref{eq:1} remains to consider. First note that \eqref{eq:1} can be written as
    \begin{align*}
        &\E_1\left(\left(f(h_n(X_n))^*\E_2 g(y_n(X_n, M_n))^*\right) -  f(h_n(X_n))_* \E_2g(y_n(X_n, M_n))_*\right)
        \\&= \E_1\left(f(h_n(X_n))^*\left(\E_2 g(y_n(X_n, M_n))^* - \E_2 g(y_n(X_n, M_n))_*\right)\right) \\&\quad + \E_1\left(\left(f(h_n(X_n))^* - f(h_n(X_n))_*\right)\E_2 g(y_n(X_n, M_n))_*\right)
        \\&\leq ||f||_{\infty}\E_1\left(\E_2 g(y_n(X_n, M_n))^* - \E_2 g(y_n(X_n, M_n))_*\right) \\&\quad + ||g||_{\infty}\E_1\left(f(h_n(X_n))^* - f(h_n(X_n))_*\right)
    \end{align*}
    By \eqref{eq:precise2} and the dominated convergence theorem, the first summand converges to zero. The second summand converges to zero since $h_n(X_n)$ is asymptotically measurable.  Consequently, $\E^* f(h_n(X_n))g(y_n(X_n, M_n)) \to \E f(H) \E g(Y) $ follows.
\end{proof}

The following theorem is an extension of Theorem~3.10.11 in \cite{vdVW23}, where $\mathbb P_n, \hat{\mathbb P}_n$ may be arbitrary maps instead of random elements, different limits $\mathbb G$ and $\widehat{\mathbb G}$ are allowed for the empirical process and its resampling counterpart, and the centering element, say $P_n$, may depend on $n$. 
This theorem is in particular applicable for $\widehat{\mathbb P}_n = \widehat{\mathbb P}_n(X_n,M_n)$ being the permutation empirical measure, i.e., $\mathbb{P}^{\pi}_{\boldsymbol{n}} := (\mathbb{P}^{\pi}_{1,n_1},\dots,\mathbb{P}^{\pi}_{m,n_m})$, or the pooled bootstrap empirical measure, i.e., $\hat{\mathbb{P}}_{\boldsymbol{n}} := (\hat{\mathbb{P}}_{1,n_1},\dots,\hat{\mathbb{P}}_{m,n_m})$. Here,
$X_n$ denotes the data and $M_n$ denotes the randomness of the resampling procedures. However, we do not restrict to the cases of permutation and pooled bootstrap empirical processes in the following theorem but allow more general processes $\widehat{\mathbb P}_n$. 

\begin{theorem}[Conditional Delta-Method]\label{DeltaMethod}
    Let $(\mathbb D, ||.||_{\mathbb D}), (\mathbb E, ||.||_{\mathbb E})$ be normed spaces, $X_n:\Omega_1\to\chi_{1n}, M_n:\Omega_2\to\chi_{2n}$ be sequences of functions, where $(\Omega_1\times\Omega_2, \mathcal{A}_1\times\mathcal{A}_2, Q_1 \otimes Q_2)$ denotes a product probability space. Furthermore, let
    $\phi:\mathbb D_{\phi} \subset \mathbb D \to \mathbb E$ be uniformly Hadamard differentiable at $P\in\mathbb D_{\phi}$ tangentially to a subspace $\mathbb D_0 \subset \mathbb D$.
    Moreover, let  $r_n$ be a sequence of constants tending to infinity,
    $P_n$ be a sequence in $\mathbb D_{\phi}$ with $P_n\to P$ and 
    $\mathbb P_n = \mathbb P_n(X_n) : \Omega_1 \to \mathbb{D}_\phi$ a map
    with $r_n (\mathbb P_n-P_n)\rightsquigarrow\mathbb G$ for some separable Borel measurable random element $\mathbb G: \Omega_1\to\mathbb D_0$. Additionally, let 
    $\widehat{\mathbb P}_n = \widehat{\mathbb P}_n(X_n,M_n):\Omega_1\times\Omega_2 \to \mathbb D_{\phi}$ a map
    with
    \begin{align}\label{eq:condConvergence}
        r_n (\widehat{\mathbb P}_n - \mathbb P_n)\rightsquigarrow \widehat{\mathbb G}
    \end{align}
    conditionally on $X_n$ in outer probability for some separable Borel measurable random element $\widehat{\mathbb G}:\Omega_1\times\Omega_2\to\mathbb D_0$. Then, we have 
    $r_n \left(\phi(\widehat{\mathbb P}_n) - \phi(\mathbb P_n)\right) \rightsquigarrow \phi_{P}^{\prime}(\widehat{\mathbb G})$ conditionally on $X_n$ in outer probability.
\end{theorem}

\begin{proof}
    We proceed analogously as in the proof of Theorem~3.10.11 in \cite{vdVW23} by applying their Theorem~3.10.5 (rather than Theorem~3.10.4).
    First note that we may assume without loss of generality that the derivative $\phi_P^{\prime}:\mathbb D \to \mathbb E$ is not only defined and continuous but also linear on the whole space $\mathbb D$ by their Problem~3.10.18. 
    For all $h\in BL_1(\mathbb E)$, we have $h \circ \phi_P^{\prime} \in BL_{\max\{1,||\phi_P^{\prime}||\}}(\mathbb D)$ where $||\phi_P^{\prime}||$ is the operator norm of $\phi_P^\prime$.
    By \eqref{eq:condConvergence}, it follows that
    \begin{align*}
        \sup\limits_{h\in BL_1(\mathbb E)} \left| \E_2 h\left(\phi_P^{\prime}\left(r_n(\widehat{\mathbb P}_n - \mathbb P_n)\right)\right)^{2*}  - \E h(\phi_P^{\prime}( \widehat{\mathbb G} ) )\right| \to 0
    \end{align*} in outer probability.
    Let $\varepsilon > 0$ be arbitrary.
    Since 
    $|h(A)^{2*} - h(B)^{2*}| \leq |h(A) - h(B)|^{*} \leq ||A - B||_{\mathbb E}^{*}$ 
    holds for all $h\in BL_1(\mathbb E), A,B \in \mathbb E$, it follows that
    \begin{align*}
        &\sup\limits_{h\in BL_1(\mathbb E)} \left| \E_2 h\left(r_n \left(\phi(\widehat{\mathbb P}_n) - \phi(\mathbb P_n)\right)\right)^{2*} - \E_2h\left(\phi_P^{\prime}\left(r_n(\widehat{\mathbb P}_n - \mathbb P_n)\right)\right)^{2*}  \right|
        \\&\leq \varepsilon + 2 Q_2\left( \left|\left|r_n \left(\phi(\widehat{\mathbb P}_n) - \phi(\mathbb P_n)\right)  - \phi_P^{\prime}\left(r_n(\widehat{\mathbb P}_n - \mathbb P_n)\right)\right|\right|_{\mathbb E}^* > \varepsilon \right).
    \end{align*}
    Theorem~3.10.5 in \cite{vdVW23} implies that
     $$\left|\left|r_n \left(\phi({\mathbb P}_n ) - \phi(P_n)\right)  - \phi_P^{\prime}\left(r_n({\mathbb P}_n - P_n)\right)\right|\right|_{\mathbb E}^* = o_{Q_1}(1).$$ 
    By Lemma~\ref{JointConvergenceLemma} and the continuous mapping theorem, we have 
    \begin{align*}
        r_n (\widehat{\mathbb P}_n - P_n) = r_n (\widehat{\mathbb P}_n - \mathbb P_n) + r_n (\mathbb P_n-P_n) \rightsquigarrow \widehat{\mathbb G} + \mathbb G
    \end{align*} unconditionally for independent $\widehat{\mathbb G}, \mathbb G$. 
    Hence, Theorem~3.10.5 in \cite{vdVW23} implies that
    $$\left|\left|r_n \left(\phi(\widehat{\mathbb P}_n ) - \phi(P_n)\right)  - \phi_P^{\prime}\left(r_n(\widehat{\mathbb P}_n - P_n)\right)\right|\right|_{\mathbb E}^* = o_{Q_1\otimes Q_2}(1)$$
    and, by the triangle inequality, that
    $$\left|\left|r_n \left(\phi(\widehat{\mathbb P}_n ) - \phi({\mathbb P}_n)\right)  - \phi_P^{\prime}\left(r_n(\widehat{\mathbb P}_n - {\mathbb P}_n)\right)\right|\right|_{\mathbb E}^* = o_{Q_1\otimes Q_2}(1).$$
    Markov's inequality yields that
    \begin{align*}
        &Q_1\left(Q_2\left( \left|\left|r_n \left(\phi(\widehat{\mathbb P}_n) - \phi(\mathbb P_n)\right)  - \phi_P^{\prime}\left(r_n(\widehat{\mathbb P}_n - \mathbb P_n)\right)\right|\right|_{\mathbb E}^* > \varepsilon \right) > \delta\right)
        \\& \leq \E_1 Q_2\left( \left|\left|r_n \left(\phi(\widehat{\mathbb P}_n) - \phi(\mathbb P_n)\right)  - \phi_P^{\prime}\left(r_n(\widehat{\mathbb P}_n - \mathbb P_n)\right)\right|\right|_{\mathbb E}^* > \varepsilon \right) / \delta
        \\& = (Q_1 \otimes Q_2) \left( \left|\left|r_n \left(\phi(\widehat{\mathbb P}_n) - \phi(\mathbb P_n)\right)  - \phi_P^{\prime}\left(r_n(\widehat{\mathbb P}_n - \mathbb P_n)\right)\right|\right|_{\mathbb E}^* > \varepsilon \right) / \delta \\&\to 0
    \end{align*}
    for all $\delta > 0$.
    Thus, it follows $$  \sup\limits_{h\in BL_1(\mathbb E)} \left| \E_2 h\left(r_n \left(\phi(\widehat{\mathbb P}_n) - \phi(\mathbb P_n)\right)\right)^{2*} - \E_2h\left(\phi_P^{\prime}\left(r_n(\widehat{\mathbb P}_n - \mathbb P_n)\right)\right)^{2*}  \right| \to 0 $$ in outer probability. Analogously, we can conclude     \begin{align}\label{eq:precise3}        \sup\limits_{h\in BL_1(\mathbb E)} \left| \E_2 h\left(r_n \left(\phi(\widehat{\mathbb P}_n) - \phi(\mathbb P_n)\right)\right)^* - \E_2h\left(\phi_P^{\prime}\left(r_n(\widehat{\mathbb P}_n - \mathbb P_n)\right)\right)^*  \right| \to 0    \end{align} in outer probability.

    For the asymptotic measurability in outer probability, write
    \begin{align*}
        &E_2h\left(r_n \left(\phi(\widehat{\mathbb P}_n) - \phi(\mathbb P_n)\right) \right)^* - E_2h\left(r_n \left(\phi(\widehat{\mathbb P}_n) - \phi(\mathbb P_n)\right) \right)_*
        \\&\leq \left| E_2h\left(r_n \left(\phi(\widehat{\mathbb P}_n) - \phi(\mathbb P_n)\right) \right)^* 
        - E_2h\left(\phi'_P\left(r_n(\widehat{\mathbb P}_n - {\mathbb P}_n)\right)\right)^*\right|
        \\&\quad+ E_2h\left(\phi'_P\left(r_n(\widehat{\mathbb P}_n - {\mathbb P}_n)\right)\right)^*
        - E_2h\left(\phi'_P\left(r_n(\widehat{\mathbb P}_n - {\mathbb P}_n)\right)\right)_*
        \\&\quad+ \left|E_2\left(-h\left(r_n \left(\phi(\widehat{\mathbb P}_n) - \phi(\mathbb P_n)\right) \right)\right)^*- E_2\left(-h\left(\phi'_P\left(r_n(\widehat{\mathbb P}_n - {\mathbb P}_n)\right)\right)\right)^*
         \right|.
    \end{align*}
    Then, it is easy to obtain the asymptotic measurability \eqref{eq:precise2} in outer probability from \eqref{eq:condConvergence} and \eqref{eq:precise3}.
\end{proof}

The following assertions will be stated in terms of the permutation and pooled bootstrap empirical measures.

\begin{cor}[Conditional Delta-Method for the Permutation Empirical Process]\label{cor1}
Let $\mathcal{F}$ satisfy $||P_i||_{\mathcal{F}} < \infty$ for all $i=1,...,m$ and $(\mathbb E, ||.||_{\mathbb E})$ be a normed space, 
    $\phi:(\ell^{\infty}(\mathcal F))^m \to \mathbb E$ be uniformly Hadamard differentiable at $(H,\dots,H)\in (\ell^{\infty}(\mathcal F))^m$ 
    tangentially to $(\ell^{\infty}(\mathcal F))^m$, where again $H = \sum_{i=1}^m \lambda_i P_i$.
    Then we have, as $\min_{j=1,\dots, m} n_j \to \infty$, 
    $$\sqrt{N} \left(\phi(\mathbb{P}_{1,n_1}^{\pi},\dots,\mathbb{P}_{m,n_m}^{\pi}) - \phi(\mathbb{H}_N,\dots,\mathbb{H}_N)\right) \rightsquigarrow \phi_{(H,\dots,H)}^{\prime}(\mathbb{G}_H^{\pi})$$ conditionally on \eqref{eq:data} in outer probability.
\end{cor}
\begin{proof}
This is an application of Theorem~\ref{DeltaMethod} together with \eqref{eq:Hn} and the conditional convergence of the permutation empirical process, see Theorem~\ref{Permutation}.
\end{proof}

The very same holds for the pooled bootstrap.

\begin{cor}[Conditional Delta-Method for the Pooled Bootstrap Empirical Process]\label{cor2}
   Let $\mathcal{F}$ satisfy $||P_i||_{\mathcal{F}} < \infty$ for all $i=1,...,m$ and $(\mathbb E, ||.||_{\mathbb E})$ be a normed space, 
    $\phi:(\ell^{\infty}(\mathcal F))^m \to \mathbb E$ be uniformly Hadamard differentiable at $(H,\dots,H)\in (\ell^{\infty}(\mathcal F))^m$ tangentially to $(\ell^{\infty}(\mathcal F))^m$.
Then we have, as $\min_{j=1,\dots, m} n_j \to \infty$, 
$$\sqrt{N} \left(\phi(\hat{\mathbb{P}}_{1,n_1},\dots,\hat{\mathbb{P}}_{m,n_m}) - \phi(\mathbb{H}_N,\dots,\mathbb{H}_N)\right) \rightsquigarrow \phi_{(H,\dots,H)}^{\prime}(\lambda_1^{-1/2}\mathbb{G}_{H,1},\dots,\lambda_m^{-1/2}\mathbb{G}_{H,m})$$ conditionally on \eqref{eq:data} in outer probability.
\end{cor}

Note that the functional of the permutation and pooled bootstrap empirical processes in Corollaries~\ref{cor1} and~\ref{cor2} cannot mimic the same limit distribution as the functional of the original empirical process in Theorem~\ref{thm2}; this is similar as for the randomization empirical process in \cite{dobler2023}.
Hence, the corollaries are not directly applicable for inference methodologies on $\phi(\mathbf{P})$ due to altered (co-)variance structures. However, a studentization can yield the consistency of the permutation and pooled bootstrap techniques with asymptotically pivotal distributions in many cases; cf.\ \cite{dobler2023} for a similar approach.

\section{Exemplary Functionals and Applications}\label{sec:Examples}
In this section, we will verify the uniform Hadamard differentiability of exemplary functionals as the Wilcoxon functional and the product integral functional.
The verifications of the uniform Hadamard differentiability of these functionals roughly follow the lines of Lemma~3.10.18 and Lemma~3.10.32 in \cite{vdVW23}.

\subsection{Wilcoxon functional}\label{ssec:Wilcoxon}
    Let $[a,b]\subset \overline{\mathbb R}$ and $$\psi: BV_M[a,b] \times BV_M[a,b] \to D[a,b], \quad \psi(A,B) := \int_{(a,.]} A\;\mathrm{d}B$$ denote the Wilcoxon functional, where $BV_M[a,b] \subset D[a,b]$ denotes the subset of càdlàg functions $[a,b] \to \mathbb{R}$ with total variation bounded by $M < \infty$ 
    equipped with the sup-norm.
We aim to show the uniform Hadamard differentiability at $(A,B)\in\mathbb{D}_{\psi}$ tangentially to $D[a,b]\times D[a,b]$ with Hadamard derivative $$ \psi'_{(A,B)}: {D}[a,b]\times {D}[a,b] \to D[a,b], \quad \psi'_{(A,B)}(\alpha,\beta) = \int_{(a,.]}A\;\mathrm{d}\beta + \int_{(a,.]}\alpha\;\mathrm{d}B ;$$ cf. Lemma~3.10.18 in \cite{vdVW23}. 
Here and below, the integral w.r.t.\ $\beta$ is defined via integration by parts if $\beta$ has unbounded variation.
Let $t\to 0$, $A_t, A, B_t,B\in BV_M[a,b]$, $\alpha_t, \alpha, \beta_t, \beta \in D[a,b]$ such that $||A_t-A||_{\infty}\to 0, ||B_t-B||_{\infty}\to 0, ||\alpha_t-\alpha||_{\infty}\to 0, ||\beta_t-\beta||_{\infty}\to 0$ and $A_t+t\alpha_t, B_t + t\beta_t \in  BV_M[a,b]$. As in the proof of Lemma~3.10.18 in \cite{vdVW23}, we consider
\begin{align*}
    &\frac{\psi(A_t+t\alpha_t,B_t+t\beta_t)-\psi(A_t,B_t)}{t}-\psi'_{(A,B)}(\alpha_t,\beta_t) \\&= \int_{(a,.]}(A_t - A)\;\mathrm{d}\beta_t + \int_{(a,.]}\alpha_t\;\mathrm{d}(B_t+t\beta_n-B).
\end{align*} The second term converges to zero by proceeding as in the proof of Lemma~3.10.18 in \cite{vdVW23}. For the first term, we apply integration by parts to obtain
\begin{align*}
   \left| \int_{(a,.]}(A_t - A)\;\mathrm{d}\beta_t \right| &=  \left|(A_t - A)(.)\beta_t(.) - (A_t - A)(a)\beta_t(a) - \int_{(a,.]}\beta_{t}(u-)\;\mathrm{d}(A_t - A)(u)\right|
   \\&\leq 2 ||A_t-A||_{\infty}||\beta_t||_{\infty} + \left|\int_{(a,.]}\beta_{t}(u-)\;\mathrm{d}(A_t - A)(u)\right|,
\end{align*} where here and throughout $\beta_t(u-) := \lim_{s \nearrow u} \beta_t(s)$ denotes the left-continuous version of $\beta_t$ at $u$. The first term converges to zero by $||A_t-A||_{\infty}\to 0$ and the second term converges to zero as in the proof of Lemma~3.10.18 in \cite{vdVW23}. Hence, we showed that
\begin{align*}
    &\frac{\psi(A_t+t\alpha_t,B_t+t\beta_t)-\psi(A_t,B_t)}{t}-\psi'_{(A,B)}(\alpha_t,\beta_t) \to 0
\end{align*} and, by the continuity of $\psi'_{(A,B)}$, the uniform Hadamard differentiability of the Wilcoxon functional follows.

\begin{example}[Wilcoxon statistic]
    Let $a = -\infty, b = \infty$. We consider the case of two independent samples $X_1,\dots,X_n \sim F$ and $Y_1,\dots,Y_m \sim G$ with empirical distribution functions $\mathbb{F}_n, \mathbb{G}_m$, respectively. The Wilcoxon statistic $\psi(\mathbb{F}_n, \mathbb{G}_m)(\infty) = \int_{(-\infty,\infty)} \mathbb{F}_n\;\mathrm{d}  \mathbb{G}_m$ is an estimator of $\psi(F, G)(\infty) = P(X_1 \leq Y_1)$. In the following, we assume $n/(n+m) \to \lambda_1 > 0, m/(n+m) \to \lambda_2 > 0$. 
    Furthermore, let us consider the $P^{X_1}$- and $P^{Y_1}$-Donsker class
    $\mathcal{F} := \{x\mapsto 1\{x\leq t\} \mid t\in\mathbb{R} \}$, cf.\ Example~2.1.3 in \cite{vdVW23}, with $||P^{X_1}||_{\mathcal{F}}, ||P^{Y_1}||_{\mathcal{F}} \leq 1$.
    As in Example~3.10.19, we get $$\sqrt{\frac{nm}{n+m}}\left(\int \mathbb{F}_n\;\mathrm{d}  \mathbb{G}_m - \int F\;\mathrm{d}  G  \right) \rightsquigarrow \sqrt{\lambda_2}\int F\;\mathrm{d}  \mathbb{G}_G + \sqrt{\lambda_1}\int \mathbb{G}_F\;\mathrm{d}  G,$$
    where $\mathbb{G}_F, \mathbb{G}_G$ denote independent tight $F$- and $G$-Brownian bridges.

    Furthermore, we get
    $$ \sqrt{\frac{nm}{n+m}} (\mathbb{H}_{n+m} - H_{n+m}) \rightsquigarrow  \sqrt{\lambda_2} \lambda_1 \mathbb{G}_F + \sqrt{\lambda_1}\lambda_2 \mathbb{G}_G \quad\text{in $D[-\infty,\infty]$}$$
    and ${H}_{n+m} \to H := \lambda_1 F + \lambda_2 G $
    for the pooled empirical distribution function $\mathbb{H}_{n+m} := \frac{n}{n+m}\mathbb{F}_n + \frac{m}{n+m}\mathbb{G}_m$ and $H_{n+m} := \frac{n}{n+m}F + \frac{m}{n+m}G$.

    For deriving the asymptotic behaviour of the permutation and pooled bootstrap counterpart of the Wilcoxon statistic, we denote the empirical distribution functions of the permutation and pooled bootstrap samples as $\mathbb{F}_n^{\pi}, \mathbb{G}_m^{\pi}, \hat{\mathbb{F}}_n, \hat{\mathbb{G}}_m$, respectively.
    Then, Theorem~\ref{Permutation} and \eqref{eq:Bootstrap} yield
    \begin{align*}
        \sqrt{n+m} (\mathbb{F}_n^{\pi} - \mathbb{H}_{n+m}, \mathbb{G}_m^{\pi} - \mathbb{H}_{n+m}) &\rightsquigarrow (\mathbb{G}^{\pi}_{H,1},\mathbb{G}^{\pi}_{H,2}) \\
        \sqrt{n+m} (\hat{\mathbb{F}}_n - \mathbb{H}_{n+m}, \hat{\mathbb{G}}_m - \mathbb{H}_{n+m}) &\rightsquigarrow (\lambda_1^{-1/2}\mathbb{G}_{H,1} , \lambda_2^{-1/2}\mathbb{G}_{H,2})
    \end{align*} in $(D[-\infty,\infty])^2$ conditionally in outer probability, where $\mathbb{G}^{\pi}_H$ denotes a tight zero-mean Gaussian process with $$\E\left[\mathbb{G}^{\pi}_{H,i}(s)\mathbb{G}^{\pi}_{H,j}(t)\right] = (\lambda_i^{-1}1\{i=j\} -1) (H(\min\{s,t\}) - H(s)H(t))$$ and $\mathbb{G}_{H,1} , \mathbb{G}_{H,2}$ denote independent tight $H$-Brownian bridges.
    Since $H,H_{n+m},\mathbb{H}_{n+m},\mathbb{F}_n^{\pi},\mathbb{G}_n^{\pi},\hat{\mathbb{F}}_n,\hat{\mathbb{G}}_n$ are (empirical) distribution functions, the total variations are bounded by $M = 1$. By Theorem~\ref{DeltaMethod} and the uniform Hadamard differentiability of the Wilcoxon functional, we get 
    \begin{align*}
        \sqrt{n+m} (\psi(\mathbb{F}_n^{\pi},\mathbb{G}_m^{\pi}) - \psi(\mathbb{H}_{n+m},  \mathbb{H}_{n+m})) &\rightsquigarrow \psi_{(H,H)}'(\mathbb{G}^{\pi}_{H,1},\mathbb{G}^{\pi}_{H,2}) \\
        \sqrt{n+m} (\psi(\hat{\mathbb{F}}_n, \hat{\mathbb{G}}_m) - \psi(\mathbb{H}_{n+m},  \mathbb{H}_{n+m})) &\rightsquigarrow \psi_{(H,H)}'(\lambda_1^{-1/2}\mathbb{G}_{H,1} , \lambda_2^{-1/2}\mathbb{G}_{H,2})
    \end{align*} in $(D[-\infty,\infty])^2$ conditionally in outer probability. Thus, it follows that
    \begin{align*}
        \sqrt{\frac{nm}{n+m}} \left(\int \mathbb{F}_n^{\pi}\;\mathrm{d}\mathbb{G}_m^{\pi} - \int\mathbb{H}_{n+m}\;\mathrm{d} \mathbb{H}_{n+m})\right) &\rightsquigarrow \sqrt{\lambda_1\lambda_2}\left(\int H \;\mathrm{d}\mathbb{G}^{\pi}_{H,2} + \int \mathbb{G}^{\pi}_{H,1} \;\mathrm{d} H \right) \\
        \sqrt{\frac{nm}{n+m}} \left(\int \hat{\mathbb{F}}_n\;\mathrm{d} \hat{\mathbb{G}}_m - \int\mathbb{H}_{n+m}\;\mathrm{d}  \mathbb{H}_{n+m}\right) &\rightsquigarrow \sqrt{\lambda_1} \int H \;\mathrm{d}\mathbb{G}_{H,2} + \sqrt{\lambda_2} \int \mathbb{G}_{H,1} \;\mathrm{d} H
    \end{align*} conditionally in outer probability by Slutsky’s lemma.
\end{example}

\begin{example}[Nelson-Aalen estimator]\label{NelsonAalen} Let us consider a survival setup with multiple samples, i.e., each data point $\mathbf{X}_{ij} = (Z_{ij},\Delta_{ij})$ consists of a (censored) failure time $Z_{ij} = \min\{X_{ij},C_{ij}\}$ and a censoring status $\Delta_{ij} = 1\{X_{ij} \leq C_{ij}\}$, see Example~3.10.20 in \cite{vdVW23} for details. Furthermore, let
\begin{align*}
    \overline{\mathbb H}_{i,n_i} (t) := \frac{1}{n_i} \sum_{j=1}^{n_i}  1\{Z_{ij}\geq t\} \quad\text{and} \quad 
    \mathbb{H}_{i,n_i}^{uc} (t) := \frac{1}{n_i} \sum_{j=1}^{n_i} \Delta_{ij} 1\{Z_{ij}\leq t\}
\end{align*} denote the survival function of the observation times and the empirical subdistribution functions of the uncensored failure times, respectively, and $\overline{H}_i(t):= P(Z_{ij}\geq t), {H}_{i}^{uc}(t):= P(Z_{ij} \leq t, \Delta_{ij} = 1)$.
The Nelson-Aalen estimator $$ \Lambda_{i,n_i}(.) := \int_{[0,.]} \frac{1}{\overline{\mathbb H}_{i,n_i}}\;\mathrm{d}\mathbb{H}_{i,n_i}^{uc}$$ estimates the cumulative hazard function $$ \Lambda_{i}(.) := \int_{[0,.]} \frac{1}{\overline{H}_{i}}\;\mathrm{d}{H}_{i}^{uc}. $$

To derive the asymptotic behaviour of the Nelson-Aalen estimators, let us consider the $P_i := P^{(Z_{i1},\Delta_{i1})}$-Donsker class
$$ \mathcal{F} := \left\{f_t^{(1)}: (x,d)\mapsto  {1}\{x \geq t\}, f_t^{(2)}: (x,d)\mapsto {1}\{x \leq t, d = 1\} \mid t\in[0,\tau]\right\}$$ for some $\tau > 0$.
Then, we have
\begin{align*}
    \sqrt{N} ( \overline{\mathbb H}_{i,n_i} - \overline{H}_i , \mathbb{H}_{i,n_i}^{uc} - {H}_{i}^{uc} )_{i=1,\dots,m} \rightsquigarrow (\lambda_i^{-1/2}\overline{\mathbb{G}}_i, \lambda_i^{-1/2}\mathbb{G}^{uc}_i)_{i=1,\dots,m} \quad \text{in $(\tilde{D}[0,\tau] \times D[0,\tau])^{m}$}
\end{align*}
as in Example~3.10.20 in \cite{vdVW23}. Here, $(\overline{\mathbb{G}}_i, \mathbb{G}^{uc}_i), i=1,\dots,m,$ denote independent tight, zero-mean Gaussian processes with covariance structure
\begin{align*}
    \E\left[\overline{\mathbb{G}}_i(s)\overline{\mathbb{G}}_i(t)\right] &= \overline{H}_i(\max\{s,t\}) - \overline{H}_i(s)\overline{H}_i(t), \\
    \E\left[\mathbb{G}^{uc}_i(s)\overline{\mathbb{G}}_i(t)\right] &=  ({H}_{i}^{uc}(s)-{H}_{i}^{uc}(t-))1\{t\leq s\} - {H}_{i}^{uc}(s)\overline{H}_i(t),\\
    \E\left[\mathbb{G}^{uc}_i(s)\mathbb{G}^{uc}_i(t)\right] &= {H}_{i}^{uc}(\min\{s,t\}) - {H}_{i}^{uc}(s){H}_{i}^{uc}(t) 
\end{align*} and throughout $\tilde{D}[0,\tau]$ denotes the subset of all functions $[0,\tau] \to \mathbb{R}$ that are everywhere left-continuous and have right limits everywhere, equipped with the sup-norm.

For the pooled empirical subdistribution functions $ \overline{\mathbb{H}}_{N} :=  \sum_{i=1}^m \frac{n_i}{N} \overline{\mathbb{H}}_{i,n_i}, \mathbb{H}_{N}^{uc} :=  \sum_{i=1}^m \frac{n_i}{N} \mathbb{H}_{i,n_i}^{uc}$, it follows that
\begin{align*}
    \sqrt{N} \left(\overline{\mathbb{H}}_{N} -  \sum_{i=1}^m \frac{n_i}{N} \overline{H}_i ,  \mathbb{H}_{N}^{uc} -  \sum_{i=1}^m \frac{n_i}{N} {H}_{i}^{uc} \right) \rightsquigarrow \left(  \sum_{i=1}^m  \lambda_i^{1/2}\overline{\mathbb{G}}_i ,  \sum_{i=1}^m\lambda_i^{1/2}\mathbb{G}^{uc}_i\right) \quad \text{in $\tilde{D}[0,\tau] \times D[0,\tau]$}.
\end{align*} Furthermore, we have that $\sum_{i=1}^m \frac{n_i}{N} \overline{H}_i \to \overline{H} := \sum_{i=1}^m \lambda_i \overline{H}_i$ and $\sum_{i=1}^m \frac{n_i}{N} {H}_{i}^{uc} \to H^{uc} := \sum_{i=1}^m \lambda_i {H}^{uc}_i$ in the sup-norm.

Let us assume $\overline{H}_i(\tau) > 0$ for all $i=1,\dots,m$ in the following.
Then, the (classical) functional delta-method (Theorem~3.10.4 in \cite{vdVW23}) implies 
\begin{align*}
    \sqrt{N} ( \Lambda_{i,n_i} - \Lambda_i  )_{i=1,\dots,m}  \rightsquigarrow \left(\lambda_i^{-1/2} \mathbb{Z}_i(C_i) \right)_{i=1,\dots,m} \quad \text{in $(D[0,\tau])^{m}$},
\end{align*}
where $\mathbb{Z}_i, i=1,\dots,m,$ are independent standard Brownian motions and $$C_i(.) = \int_{[0,.]} \frac{1 - \Delta\Lambda_i}{\overline{H}_i}\;\mathrm{d}\Lambda_i.$$
 as in Example~3.10.20 in \cite{vdVW23}.

Now, we are considering the permutation and pooled bootstrap counterparts of the Nelson-Aalen estimators.
Therefore, we denote all processes and statistics introduced above with a $\pi$ in the superscript, if they are based on the permuted data $\boldsymbol{Z}_{NR_1},\dots,\boldsymbol{Z}_{NR_N}$ instead of the original data, and with a hat $\hat{\phantom{x}}$, if they are based on the bootstrapped data $\hat{\boldsymbol{Z}}_{N1},\dots,\hat{\boldsymbol{Z}}_{NN}$.
Theorem~\ref{Permutation} and \eqref{eq:Bootstrap} imply the conditional weak convergence of the permutation and pooled bootstrap empirical processes \begin{align*}
    \sqrt{N} (  \overline{\mathbb H}_{i,n_i}^{\pi} - \overline{\mathbb H}_{N}, \mathbb{H}_{i,n_i}^{uc,\pi} - \mathbb{H}_{N}^{uc} )_{i=1,\dots,m} &\rightsquigarrow ( \overline{\mathbb{G}}_i^{\pi}, \mathbb{G}^{uc,\pi}_i)_{i=1,\dots,m} \quad \text{in $(\tilde{D}[0,\tau] \times D[0,\tau])^{m}$} ,\\
    \sqrt{N} (  \hat{\overline{\mathbb H}}_{i,n_i} -  \overline{\mathbb H}_{N} , \hat{\mathbb{H}}_{i,n_i}^{uc} - {\mathbb{H}}_{N}^{uc})_{i=1,\dots,m} &\rightsquigarrow (\lambda_i^{-1/2}\hat{\overline{\mathbb{G}}}_i, \lambda_i^{-1/2}\hat{\mathbb{G}}^{uc}_i)_{i=1,\dots,m} \quad \text{in $(\tilde{D}[0,\tau] \times D[0,\tau])^{m}$},
\end{align*}
conditionally in outer probability due to $||P_i||_{\mathcal F} \leq 1$. Here, $(\overline{\mathbb{G}}_i^{\pi}, \mathbb{G}^{uc,\pi}_i)_{i=1,\dots,m}$ is a tight, zero-mean Gaussian process with
\begin{align*}
    \E\left[{\overline{\mathbb{G}}}_i^{\pi}(s){\overline{\mathbb{G}}}^{\pi}_j(t)\right] &= (\lambda_i^{-1}1\{i=j\} -1)(\overline{H}(\max\{s,t\}) - \overline{H}(s)\overline{H}(t) )\\
    \E\left[\mathbb{G}^{uc,\pi}_i(s){\overline{\mathbb{G}}}_j^{\pi}(t)\right] &=  (\lambda_i^{-1}1\{i=j\} -1)\left(({H}^{uc}(s)-{H}^{uc}(t-))1\{t\leq s\} - {H}^{uc}(s)\overline{H}(t)\right) \\
    \E\left[{\mathbb{G}}_i^{uc,\pi}(s)\mathbb{G}^{uc,\pi}_j(t)\right] &= (\lambda_i^{-1}1\{i=j\} -1)({H}^{uc}(\min\{s,t\}) - {H}^{uc}(s){H}^{uc}(t))
\end{align*}
and $( \hat{\overline{\mathbb{G}}}_i, \hat{\mathbb{G}}^{uc}_i), {i=1,\dots,m},$ are independent tight, zero-mean Gaussian processes with
\begin{align*}
    \E\left[\hat{\overline{\mathbb{G}}}_i(s)\hat{\overline{\mathbb{G}}}_i(t)\right] &= \overline{H}(\max\{s,t\}) - \overline{H}(s)\overline{H}(t) \\
    \E\left[\hat{\mathbb{G}}^{uc}_i(s)\hat{\overline{\mathbb{G}}}_i(t)\right] &=  ({H}^{uc}(s)-{H}^{uc}(t-))1\{t\leq s\} - {H}^{uc}(s)\overline{H}(t)\\
    \E\left[\hat{\mathbb{G}}_i^{uc}(s)\hat{\mathbb{G}}^{uc}_i(t)\right] &= {H}^{uc}(\min\{s,t\}) - {H}^{uc}(s){H}^{uc}(t) .
\end{align*}

The Nelson-Aalen functional is a composition of the functionals
$$\tilde{\psi}: \tilde{BV}_M[0,\tau] \times BV_M[0,\tau] \to D[0,\tau], \quad \tilde{\psi} (A,B) := \int_{[0,.]}A\;\mathrm{d}B$$ and $(A,B) \mapsto (1/A,B)$, where here and throughout $\tilde{BV}_M[0,\tau] \subset \tilde{D}[0,\tau]$ denotes the subset of functions with total variation bounded by $M < \infty$. 
Furthermore, we set $\Delta B (0) := B (0)$ for $B\in D[0,\tau]$ in the following to guarantee a well-defined jump in $0$.
Similarly to the above calculations for the Wilcoxon functional, it can be shown that $\tilde{\psi}$ is uniformly Hadamard differentiable at $(A,B) \in \mathbb{D}_{\tilde{\psi}}$ with Hadamard derivative
$$ \tilde{\psi}'_{(A,B)}: \tilde{D}[0,\tau]\times {D}[0,\tau] \to D[0,\tau], \quad \tilde{\psi}'_{(A,B)}(\alpha,\beta) = \int_{[0,.]}A\;\mathrm{d}\beta + \int_{[0,.]}\alpha\;\mathrm{d}B .$$ Furthermore, it is easy to show that $(A,B) \mapsto (1/A,B)$ is uniformly Hadamard differentiable at $(A,B)\in \tilde{D}[0,\tau]\times {D}[0,\tau]$ such that $|A| \geq \varepsilon$ for some $\varepsilon > 0$
with Hadamard derivative $(\alpha,\beta) \mapsto (-\alpha/A^2, \beta)$.
Hence, the Nelson-Aalen functional $(A,B) \mapsto \tilde{\psi}(1/A, B)$ is uniformly Hadamard differentiable at $(A,B)$ with Hadamard derivative $(\alpha,\beta) \mapsto \tilde{\psi}'_{(1/A,B)}(-\alpha/A^2, \beta)$ by the chain rule (Theorem~\ref{chain}), where $|A| \geq \varepsilon$ and $(1/A,B)\in\mathbb{D}_{\tilde{\psi}}$.

Since $\overline{H}_i, H_i^{uc}, \overline{H}, H^{uc}$ are positive monotone functions, we have $$\overline{H}_i, \overline{H} > \min\{\overline{H}_1(\tau),\dots, \overline{H}_m(\tau)\} =: 2\varepsilon > 0$$ and the total variation of $1/\overline{H}_i, H_i^{uc}, 1/\overline{H}, H^{uc}$ is bounded by $M := \varepsilon^{-1}$. 
Moreover, the considered empirical processes are contained in $\{A \mid A \geq \varepsilon, 1/A \in \tilde{BV}_M[0,\tau]\} \times BV_M[0,\tau]$ with probability tending to 1 by monotonicity and Glivenko-Cantelli arguments. 
Hence, the uniform Hadamard differentiability of the Nelson-Aalen functional and the conditional delta-method (Theorem~\ref{DeltaMethod}) yield 
\begin{align*}
    \sqrt{N} ( \Lambda_{i,n_i}^{\pi} -  \Lambda_N )_{i=1,\dots,m}  &\rightsquigarrow \left(\mathbb{Z}_i^{\pi} \right)_{i=1,\dots,m} \quad \text{in $(D[0,\tau])^{m}$}, \\
    \sqrt{N} ( \hat{\Lambda}_{i,n_i} - \Lambda_N  )_{i=1,\dots,m}  &\rightsquigarrow \left(\lambda_i^{-1/2} \mathbb{Z}_i(C) \right)_{i=1,\dots,m} \quad \text{in $(D[0,\tau])^{m}$}
\end{align*}
similarly to the calculations in Example~3.10.20 in \cite{vdVW23}, where 
$\Lambda_N(.) := \int_{[0,.]} \frac{1}{\overline{\mathbb H}_{N}}\;\mathrm{d}\mathbb{H}_{N}^{uc}$ denotes the pooled Nelson-Aalen estimator, $$ C(.) := \int_{[0,.]} \frac{1 - \Delta\Lambda}{\overline{H}}\;\mathrm{d}\Lambda,$$ $\Lambda(.) := \int_{[0,.]} \frac{1}{\overline{H}}\;\mathrm{d}{H}^{uc}$ 
and $\left(\mathbb{Z}_i^{\pi} \right)_{i=1,\dots,m}$ is a zero-mean Gaussian process with 
\begin{align*}
    \E\left[{\mathbb{Z}^{\pi}}_i(s){{\mathbb{Z}}}_j^{\pi}(t)\right] &=  (\lambda_i^{-1}1\{i=j\} -1) C(\min\{s,t\}).
\end{align*}
\end{example}

\subsection{Product integral}\label{ssec:Product}
Consider $\phi$, the product-integral functional, i.e., 
$$\phi: BV_{M}^{>-1}[a,b] \to D[a,b], \quad A \mapsto \Prodi_{u \in (a,\cdot]}(1 + A(du)). $$
Here, $BV_{M}^{>-1}[a,b] \subset D[a,b]$ is the subset of functions $[a,b] \to \mathbb{R}$ with total variation bounded by $M$ and whose jumps are contained in $(-1, \infty)$ and bounded away from $-1$.
To analyze the uniform Hadamard differentiability of $\phi$, let $t_n \to 0$,  $A_n , A \in BV_{M}^{>-1}[a,b]$ such that $\|A_n - A \|_\infty\to 0$, and $\alpha_n ,  \alpha \in D[a,b] $ such that $\|\alpha_n - \alpha\|_\infty\to 0$ and $A_n + t_n \alpha_n \in BV_{M}^{>-1}[a,b]$.
Let $\varepsilon > 0$ and $\tilde \alpha \in BV[a,b]$ such that $\|\alpha - \tilde \alpha \|_\infty < \varepsilon$, $\| \alpha_n - \alpha\|_\infty < \varepsilon$, and $\| A_n - A \|_\infty < \varepsilon$ for sufficiently large $n$.
This function $\tilde \alpha$ can be defined piece-wise constant and with finitely many jumps because it approximates the c\`adl\`ag function $\alpha$. Also, because the sequence $(\alpha_n)_n$ approximates $\alpha$ uniformly, it is clear that such a function $\tilde \alpha$ exists.

It is well-known that 
\begin{align*}
\phi'_A : D[a,b] \to D[a,b] , \alpha \mapsto \int_{(a,\cdot]} \phi(A)(u-) \alpha(du) \frac{\phi(A)(\cdot)}{\phi(A)(u)} =  \int_{(a,\cdot]} \frac{\alpha(du)}{1+\Delta A(u)}\phi(A)(\cdot)
\end{align*}
defines the Hadamard derivative of $\phi$ at $A$ in the classical sense; cf.\ \cite{gill90}, where here and throughout $\Delta A(u) := A(u) - A(u-)$ denotes the increment of $A$ at $u$.
Note that the Hadamard derivative above may also be written as
\begin{align*}
   \phi_A'(\alpha) 
& = \phi(A)(\cdot)\left(\alpha(\cdot) - \alpha(a) - \sum_{u \in D_A \cap (a, \cdot]} \frac{\Delta A(u) \Delta \alpha(u)}{1 + \Delta A(u)} \right)
\end{align*}
where $D_A \subset (a,b]$ is the set of discontinuities of $A$.
Due to the finite variation of $A$ and its boundedness of its jumps away from $-1$, this representation reveals that $\tilde \phi' := ((A,\alpha) \mapsto \phi'_A(\alpha))$ defines a continuous functional from $BV_{M}^{>-1}[a,b] \times D[a,b] \to D[a,b] $ with respect to the maximum-supremum norm.
To see this, let us focus on the sum-term and notice that
\begin{align*}
&\sum_{u \in D_{A_n} \cap (a, \cdot]} \frac{\Delta A_n(u) \Delta \alpha_n(u)}{1+\Delta A_n(u)} - \sum_{u \in D_{A} \cap (a, \cdot]} \frac{\Delta A(u) \Delta \alpha(u)}{1+\Delta A(u)} 
\\
&= \sum_{u \in (D_{A_n} \cup D_A) \cap (a, \cdot]} \frac{\Delta A_n(u) \Delta A(u) \Delta (\alpha_n(u)-\alpha(u)) + \Delta A_n(u) \Delta \alpha_n(u) - \Delta A(u) \Delta \alpha (u)}{(1+\Delta A_n(u))(1+\Delta A(u))} 
\end{align*}
Choose $\delta > 0$ sufficiently small, i.e., $\delta < \min(\varepsilon/2,\inf_u (1+\Delta A(u)))$
such that a finite, positive constant $K \geq \sup_{u}(1+\Delta A (u) - \delta)^{-1} $ exists.
Now, choose $n_0$ sufficiently large such that  $\sup_u|\Delta A_n (u) - \Delta A(u)| \leq 2 \| A_n - A \|_\infty \leq 2 \varepsilon < \delta$ for all $n \geq n_0$. 
Let us only consider such $n$ henceforth.
This implies that $K \geq \sup_{u}(1+\Delta A_n (u))^{-1} $.
Since the jumps of $A$ and $A_n$ are bounded away from $-1$,
the supremum norm of the previous display is bounded above by $K^2$ (due to the denominator) times the sum of
\begin{align*}
    & 2\|\alpha_n - \alpha \|_\infty \|A_n \|_{BV} \|A \|_{BV} && \text{ for the term $\Delta A_n(u) \Delta A(u) \Delta (\alpha_n(u)-\alpha(u))$,} \\
    & 2 \|\alpha_n - \alpha \|_\infty \|A_n \|_{BV}  && \text{ for a term  } \Delta A_n(u) \Delta ( \alpha_n(u) - \alpha(u)), \\
    & 2 \| \alpha - \tilde \alpha\|_\infty (\|A_n\|_{BV} + \| A\|_{BV})  && \text{ for a term } \Delta (A_n(u) - A(u) )\Delta ( \alpha(u) - \tilde \alpha(u)),\\
    & 2 \|A_n - A\|_\infty \| \tilde \alpha \|_{BV}   && \text{ for a term } \Delta (A_n(u) - A(u) ) \Delta\tilde \alpha(u),
\end{align*}
where $||.||_{BV}$ denotes the total variation.
Hence, an upper bound is given by $K^2 (8 \max({M}^2,1) + 2\|\tilde \alpha \|_{BV}) \varepsilon$ which is arbitrarily small because $\varepsilon > 0$ was arbitrary.

To verify the uniform Hadamard differentiability, we need to show that the following term converges to zero:
\begin{align*}
    t_n^{-1}(\phi(A_n + t_n \alpha_n) - \phi(A_n)) - \phi'_A(\alpha) =: I + II
\end{align*}
where 
\begin{align*}
    I = & \ t_n^{-1}(\phi(A_n + t_n \alpha_n) - \phi(A_n)) - \phi'_{A_n}(\alpha_n), 
\\
    II = & \ \phi'_{A_n}(\alpha_n)  - \phi'_{A}(\alpha) 
\end{align*}
The second term can be written as $\tilde \phi'(A_n,\alpha_n) - \tilde \phi' (A,\alpha)$ which goes to zero as argued above.
Hence, we focus on the first term which, by Duhamel's equation, equals:
\begin{align*}
    & t_n^{-1} \int_{(a,\cdot]} \phi(A_n + t_n \alpha_n)(u-) (A_n + t_n \alpha_n - A_n)(du) \frac{\phi(A_n)(\cdot)}{\phi(A_n)(u)} - \phi'_{A_n}(\alpha_n) \\
    & = \int_{(a,\cdot]} [\phi(A_n + t_n \alpha_n)(u-) - \phi(A_n)(u-) ] (\alpha_n - \tilde \alpha + \tilde \alpha)(du) \frac{\phi(A_n)(\cdot)}{\phi(A_n)(u)} 
\end{align*}
The part with $\alpha_n - \tilde \alpha$ is arbitrarily small, which follows from integration-by-parts, combined with the fact that the involved product-integrals have a finite variation; cf.\ the proof of Theorem~7 in \cite{gill90} for similar arguments.
The upper bound for the variation norm of the involved product-integrals can be chosen independent of $n$.
The remaining part with $\tilde \alpha$ converges to zero in supremum norm due to the uniform continuity of the product-integral functional (Theorem 7 in \citealp{gill90}) combined with $A_n + t_n \alpha_n \to A,  A_n \to A$. Indeed, combine  $\|\tilde \alpha\|_{BV} < \infty$ with $\| \phi(A_n + t_n \alpha_n)(\cdot) - \phi(A_n)(\cdot)  \|_\infty \to 0 $ and 
$\sup_{a \leq u \leq t \leq b}|  {\phi(A_n)(t)}/{\phi(A_n)(u)}| < \tilde K$ for some finite constant $\tilde K$ independent of $n$; cf.\ the inequality in (20) of \cite{gill90}.

\begin{example}[Kaplan-Meier estimator]\label{kaplan} Let us consider the setup of Example~\ref{NelsonAalen}. 
The Kaplan-Meier estimator $$ \widehat{S}_{i,n_i}(\cdot) :=  \Prodi_{u \in [0,\cdot]}(1 - \widehat{\Lambda}_{i,n_i}(du))$$ estimates the survival function $$S_{i}(\cdot) := P(X_{ij} > \cdot ) = \Prodi_{u \in [0,\cdot]}(1 - {\Lambda}_{i}(du)). $$
The (classical) functional delta-method (Theorem~3.10.4 in \cite{vdVW23}) implies 
\begin{align*}
    \sqrt{N} ( \widehat{S}_{i,n_i} - S_i  )_{i=1,\dots,m}  \rightsquigarrow \left(\lambda_i^{-1/2} \mathbb{U}_i \right)_{i=1,\dots,m} \quad \text{in $(D[0,\tau])^{m}$},
\end{align*}
where $\mathbb{U}_i, i=1,\dots,m,$ are independent 
zero-mean Gaussian processes with covariance structure
\begin{align*}
    \E\left[ \mathbb{U}_i(s)\mathbb{U}_i(t) \right] = S_i(s)S_i(t) \int_{[0,\min\{s,t\}]} \frac{1}{(1 - \Delta\Lambda_i)\overline{H}_i}\;\mathrm{d}\Lambda_i,
\end{align*}
which can be shown as in Example~3.10.33 in \cite{vdVW23}.

The Kaplan-Meier functional is a composition of the functionals $$\tilde{\phi}: BV^{>-1}_{M}[0, \tau] \to D[0,\tau], \quad A \mapsto \Prodi_{u \in [0,\cdot]}(1 + A(du)), $$  $A \mapsto -A$, and the Nelson-Aalen functional. Again, we set $\Delta A (0) := A (0)$ for $A\in D[0,\tau]$ to guarantee a well-defined jump in $0$.
Similarly as above, $\tilde{\phi}$ is uniformly Hadamard differentiable at each $A \in \mathbb{D}_{\tilde\phi}$ with Hadamard derivative
\begin{align*}
{\tilde\phi}'_A : D[0,\tau] \to D[0,\tau] , \quad \alpha \mapsto \int_{[0,\cdot]} {\tilde\phi}(A)(u-) \alpha(du) \frac{{\tilde\phi}(A)(\cdot)}{\tilde{\phi}(A)(u)} .
\end{align*} Moreover, $A \mapsto -A$ is uniformly Hadamard differentiable at each $A \in D[0,\tau]$ with Hadamard derivative $\alpha \mapsto -\alpha$  due to its linearity. Thus, the Kaplan-Meier functional $A \mapsto \tilde{\phi}(-A) $ is uniformly Hadamard differentiable at each $A$ such that $-A \in \mathbb{D}_{\tilde\phi}$ with Hadamard derivative
$\alpha \mapsto - {\tilde\phi}'_A(\alpha)$ by Theorem~\ref{chain}.
To apply the conditional delta-method in Theorem~\ref{DeltaMethod}, we need to ensure that $-\Lambda_{i}, -\Lambda$ are elements in $\mathbb{D}_{\tilde\phi}$. This can be guaranteed by assuming $S_i(\tau) = P(X_{ij} > \tau) > 0$ in the following. 
Moreover, the Nelson-Aalen estimator and its permutation and pooled bootstrap counterpart are contained in $\{A \mid -A \in \mathbb{D}_{\tilde\phi}\}$ with probability tending to 1 by monotonicity and Glivenko-Cantelli arguments. By the uniform Hadamard differentiability of the Kaplan-Meier functional, Theorem~\ref{DeltaMethod} yields 
\begin{align}\label{eq:perm}
    \sqrt{N} ( \widehat{S}_{i,n_i}^{\pi} -  \widehat{S}_{N} )_{i=1,\dots,m}  &\rightsquigarrow \left(\mathbb{U}_i^{\pi} \right)_{i=1,\dots,m} \quad \text{in $(D[0,\tau])^{m}$}, \\ \label{eq:boot}
    \sqrt{N} ( \hat{\widehat{S}}_{i,n_i} -  \widehat{S}_{N}  )_{i=1,\dots,m}  &\rightsquigarrow \left(\lambda_i^{-1/2} \hat{\mathbb{U}}_i \right)_{i=1,\dots,m} \quad \text{in $(D[0,\tau])^{m}$}
\end{align}
similarly to the calculations in Example~3.10.33 in \cite{vdVW23}, where 
$\widehat{S}_{N} := \tilde{\phi}(-\widehat{\Lambda}_N)$ denotes the pooled Kaplan-Meier estimator.
Here, $\left(\mathbb{U}_i^{\pi} \right)_{i=1,\dots,m}$ is a zero-mean Gaussian process with 
\begin{align*}
 \E\left[\mathbb{U}^{\pi}_i(s){{\mathbb{U}}}_j^{\pi}(t)\right] &=  (\lambda_i^{-1}1\{i=j\} -1) S(s)S(t) \int_{[0,\min\{s,t\}]}  \frac{1}{(1 - \Delta\Lambda)\overline{H}}\;\mathrm{d}\Lambda
\end{align*} for $S := \tilde{\phi}(-\Lambda)$ and $\hat{\mathbb{U}}_i, i=1,\dots,m,$ are independent zero-mean Gaussian processes with 
\begin{align*}
 \E\left[\hat{\mathbb{U}}_i(s){\hat{\mathbb{U}}}_j(t)\right] &=   S(s)S(t) \int_{[0,\min\{s,t\}]}  \frac{1}{(1 - \Delta\Lambda)\overline{H}}\;\mathrm{d}\Lambda.
\end{align*}
\end{example}

From this example, we can deduce Theorem~4 and Theorem~5 in the supplement of \cite{doblerpauly2018bootstrap} under $S_i(\tau) > 0, i=1,2$. The uniform Hadamard differentiability of the Wilcoxon functional completes the proofs of the consistency for the permutation and pooled bootstrap counterpart of the Mann–Whitney statistic (Theorem~2 and Theorem~3 in \citealp{doblerpauly2018bootstrap}).

Furthermore, other works on resampling in survival analysis are based on Theorem~4 and Theorem~5 in the supplement of \cite{doblerpauly2018bootstrap}. This includes the resampling tests for the restricted mean survival times (RMSTs) of \cite{ditzhaus2023} and \cite{munko2024}. For the RMSTs, the application of the continuous mapping theorem with continuous function $$(D[0,\tau))^m \ni (A_1,\dots,A_m) \mapsto \left(\int_0^{\tau} A_1(t)\;\mathrm{d}t, \dots, \int_0^{\tau} A_m(t)\;\mathrm{d}t\right)$$ only requires \eqref{eq:perm} and \eqref{eq:boot} in $(D[0,\tau))^m$ instead of $(D[0,\tau])^m$.
By replacing all intervals $[0,\tau]$ in Example~\ref{kaplan} by $[0,\tau)$, it is easy to see that the assumption $S_i(\tau-) > 0, i=1,\dots,m,$ (instead of $S_i(\tau) > 0, i=1,\dots,m$) is sufficient to ensure \eqref{eq:perm} and \eqref{eq:boot} in $(D[0,\tau))^m$. Hence, the weaker assumption $S_i(\tau-) > 0, i=1,\dots,m,$ is enough for getting the consistency for the permutation and pooled bootstrap counterparts of the RMSTs as in \cite{ditzhaus2023} and \cite{munko2024}.


\subsection{Inverse map: counterexample and additional requirements}\label{ssec:Inverse}
Let $p \in \mathbb{R}$ and $A \in D[a,b]$ nondecreasing with $A(y-) \leq p \leq A(y)$ for some $y\in (a,b]$.
Then, the inverse map $\Phi_p$ at $A$ satisfies
$$ A(\Phi_p(A)-) \leq p \leq A(\Phi_p(A)), $$ 
where the exact value of $\Phi_p(A)$ is irrelevant if there is more than one solution.
Let $\mathbb{D}_{\Phi_p}$ denote the set of all nondecreasing functions $A$ with $A(y-) \leq p \leq A(y)$ for some $y\in (a,b]$.
As shown in Lemma~3.10.21 of \cite{vdVW23}, $\Phi_p : \mathbb{D}_{\Phi_p} \to (a,b]$
is Hadamard differentiable at a function $A \in \mathbb{D}_{\Phi_p}$ that is differentiable at $\Phi_p(A) =: \xi_p \in (a,b)$ such that $A(\xi_p) = p$ with positive derivative $A'(\xi_p) > 0$, tangentially to the set of functions $\alpha\in D[a,b]$ that are continuous at $\xi_p$, with Hadamard derivative $\Phi'_{p,A}(\alpha) = -\alpha(\Phi_p(A))/A'(\Phi_p(A))$ at $\alpha$. 
However, the uniform Hadamard differentiability of the inverse map does not hold under these assumptions.

\begin{figure}[bth]
        \centering
        \includegraphics[width=0.7\linewidth]{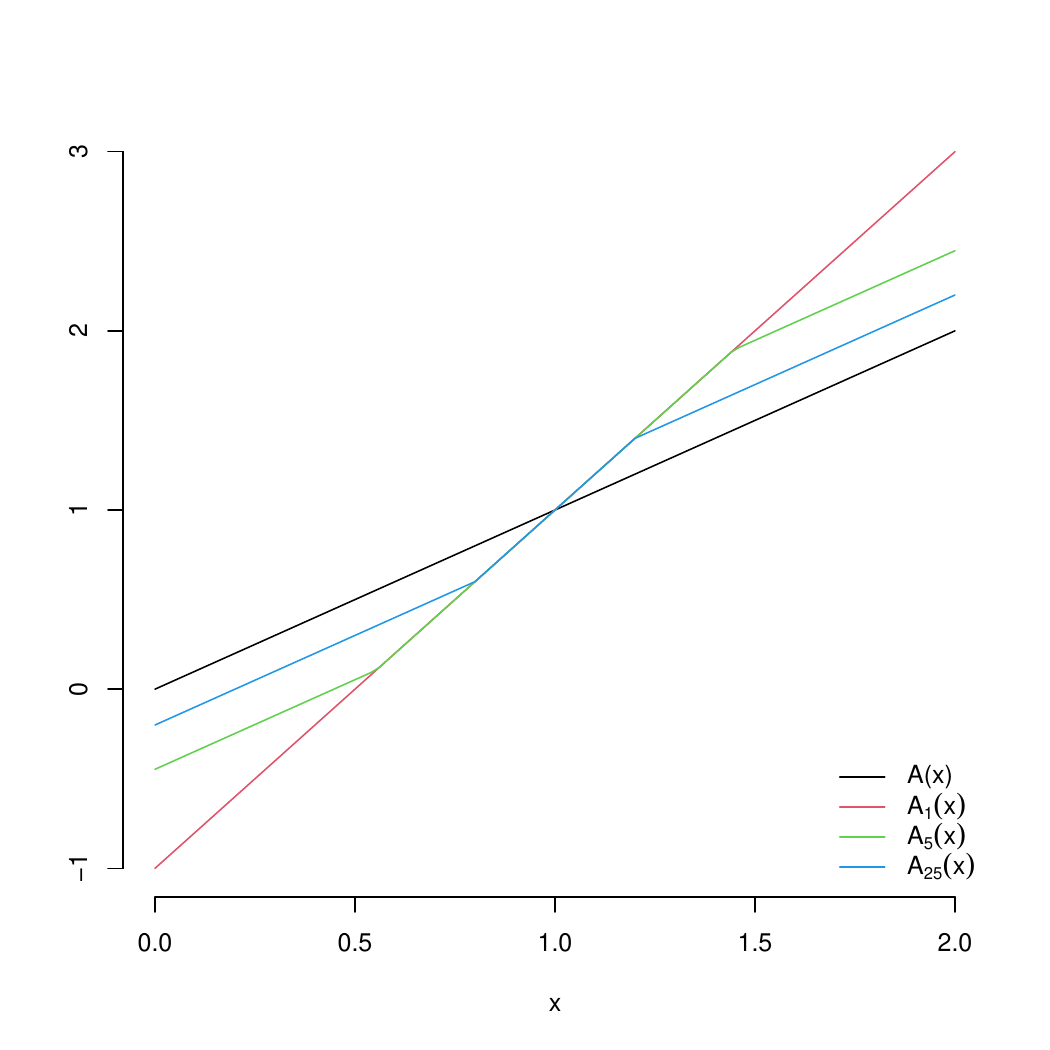}
        \caption{Illustration of the functions $A, A_1, A_5,$ and $A_{25}$.}
        \label{fig:plot}
    \end{figure}

    \begin{example}\label{counterex}
        Let $A: [0,2]\to\mathbb{R}, A(x) = x,$ \begin{align} A_n:[0,2]\to\mathbb{R},\quad  A_n(x) := 
            \begin{cases}
                x - 1/\sqrt{n} & \text{if $x \leq 1-1/\sqrt{n}$},\\
                2x - 1 & \text{if $1+1/\sqrt{n} > x > 1-1/\sqrt{n}$},\\
                x + 1/\sqrt{n} & \text{if $x \geq 1+1/\sqrt{n}$},\\
            \end{cases}
        \end{align}
$\alpha \equiv \alpha_n \equiv 1$ and $p=1$. 
An exemplary illustration for the functions can be found in Figure~\ref{fig:plot}.
For $t_n = 1/\sqrt{n}$, we have
$\Phi_p(A_n+t_n\alpha_n) = 1 - t_n/2$ since $$(A_n+t_n\alpha_n)(1-t_n/2) = A_n(1-1/(2\sqrt{n})) + 1/\sqrt{n} = 1-1/\sqrt{n} +1/\sqrt{n} = 1$$ and $\Phi_p(A_n) = 1$. Hence, $(\Phi_p(A_n+t_n\alpha_n) - \Phi_p(A_n))/t_n = -1/2$. However, $\Phi'_{p,A}(\alpha) = -\alpha(\Phi_p(A))/A'(\Phi_p(A)) = -1 \neq -1/2$ and, thus, $\Phi_p:\mathbb{D}_{\Phi_p} \to (0,2]$ is not uniformly Hadamard differentiable at $A$ tangentially to $\alpha \equiv 1$.   
    \end{example}

A restriction of the definition space $\mathbb{D}_{\Phi_p}$ can lead to uniform Hadamard differentiability. In view of Example~\ref{counterex}, such a restriction needs to exclude the rather simple sequence $(A_n)_n$. That is why such a restriction is not really of interest for applications.

However, in Lemma~S.5 in \cite{ditzhaus2021qanova}, a version of uniform Hadamard differentiability of the inverse map is shown under stricter conditions, where the rate of convergence of the converging sequence $A_n \to A$ and its increments around $\xi_p$ also needs to be controlled.
Note that for Example~\ref{counterex}, the condition \eqref{eq:incrementcondition} does not hold.

\begin{lemma}[Lemma~S.5 of \cite{ditzhaus2021qanova}]
    Let $A_n, A \in \mathbb{D}_{\Phi_p}$ such that $A$ is continuously differentiable at $\xi_p \in\mathbb R$ with positive
derivative $A'(\xi_p) > 0$. Suppose that $ \sqrt{n} ||A_n - A||_{\infty} \leq M $ for some $M > 0$
and 
\begin{align}\label{eq:incrementcondition}
    \sqrt{n} \sup_{|x|\leq K/\sqrt{n}} | A_n(\xi_p+x)-A_n(\xi_p) -A(\xi_p+x)+A(\xi_p)  | \to 0
\end{align}
for every $K > 0$. Then,
$$ \sqrt{n} \left(\Phi_p(A_n + n^{-1/2} \alpha_n) - \Phi_p(A_n)  \right) \to \Phi'_{p,A}(\alpha), $$
 where $A_n + n^{-1/2} \alpha_n \in \mathbb{D}_{\Phi_p} $ and $\alpha_n \to \alpha$ such that $\alpha$ is bounded and continuous at $\xi_p$. 
\end{lemma}

As shown in Section~S3.3.1 in the supplement of \cite{ditzhaus2021qanova}, the required conditions are fulfilled for applications on empirical distribution functions. Hence, a central limit theorem for permutation quantiles follows as shown in Lemma~S.1 in the supplement of \cite{ditzhaus2021qanova}. Analogously, a central limit theorem for pooled bootstrap quantiles could be followed, where $\hat{\gamma}(c,d) = \lambda_c^{-1/2} 1\{c = d\}$ replaces $\gamma^{\pi}(c,d) = \lambda_c^{-1} 1\{c=d\} -1 $ in the covariance formula given in Lemma~S.1.


\section{Discussion}
\label{sec:Discussion}
In this paper, we closed the gap of a suitable delta-method for resampling procedures in multiple sample problems as, e.g., the permutation and pooled bootstrap. 
Therewith, the weak convergence of uniform Hadamard differentiable functionals of resampling empirical processes can be derived in outer probability which is sufficient for most statistical applications.
In view of the extensions of the classical delta-method for, e.g., quasi-Hadamard differentiable functionals \citep{BEUTNER20102452, BeutnerZaehleBootstrap, beutner2016functionaldeltamethodbootstrapuniformly} and directionally differentiable functionals \citep{fang2019inference}, it could be studied in future research whether conditional delta-methods for more general functionals can be achieved that are applicable for resampling procedures in multiple sample problems. 
Additionally, the weak convergence could also be investigated in the outer almost sure case; cf. the supplement of \citet{beutner2016functionaldeltamethodbootstrapuniformly} for an extension of the conditional delta-method outer almost surely under measurability assumptions. 
